\documentclass[11pt]{article}

\usepackage{geometry}
\usepackage[english]{babel}
\usepackage[babel]{microtype}
\usepackage{csquotes}
\usepackage{setspace}
\usepackage{fancyhdr}
\usepackage[title]{appendix}

\usepackage{amsthm, amsmath, amssymb, amsfonts, bm, bbm, mathabx}
\usepackage{comment}
\usepackage{aliascnt}

\newtheorem{theorem}{Theorem}[section]

\newaliascnt{prop}{theorem}
\newtheorem{prop}[prop]{Proposition}
\aliascntresetthe{prop}

\newaliascnt{lemma}{theorem}
\newtheorem{lemma}[lemma]{Lemma}
\aliascntresetthe{lemma}

\newaliascnt{cor}{theorem}
\newtheorem{cor}[cor]{Corollary}
\aliascntresetthe{cor}

\newaliascnt{claim}{theorem}
\newtheorem{claim}[claim]{Claim}
\aliascntresetthe{claim}

\newaliascnt{fact}{theorem}

\aliascntresetthe{fact}

\newaliascnt{conj}{theorem}
\newtheorem{conj}[conj]{Conjecture}
\aliascntresetthe{conj}

\theoremstyle{definition}

\newaliascnt{defn}{theorem}

\aliascntresetthe{defn}
\newtheorem*{defn-non}{Definition}

\newaliascnt{exam}{theorem}

\aliascntresetthe{exam}

\newaliascnt{setup}{theorem}

\aliascntresetthe{setup}

\newaliascnt{exer}{theorem}

\aliascntresetthe{exer}

\newaliascnt{ques}{theorem}

\aliascntresetthe{ques}

\newaliascnt{rmk}{theorem}

\aliascntresetthe{rmk}

\usepackage[shortlabels]{enumitem}
\newlist{Case}{enumerate}{3}
\setlist[Case, 1]{%
  label = {\bfseries Case \arabic*.},
  labelindent=1em, labelwidth=1cm, labelsep*=1em, leftmargin=!
}
\setlist[Case, 2]{%
  label = {\bfseries Subcase \arabic{Casei}.\arabic*.},
  labelindent=-1em, labelwidth=1cm, labelsep*=1em, leftmargin=!
}
\setlist[Case, 3]{%
  label = {\bfseries Subsubcase \arabic{Casei}.\arabic{Caseii}.\arabic*.},
  labelindent=-1em, labelwidth=1cm, labelsep*=1em, leftmargin=!
}

\usepackage{url}
\usepackage{booktabs}
\usepackage{graphicx}
\usepackage{float}
\usepackage{tikz}

\usepackage[colorlinks=true,linkcolor=blue,citecolor=blue,urlcolor=blue]{hyperref}
\usepackage[capitalise]{cleveref}

\crefname{theorem}{Theorem}{Theorems}
\Crefname{theorem}{Theorem}{Theorems}
\crefname{prop}{Proposition}{Propositions}
\Crefname{prop}{Proposition}{Propositions}
\crefname{lemma}{Lemma}{Lemmas}
\Crefname{lemma}{Lemma}{Lemmas}
\crefname{cor}{Corollary}{Corollaries}
\Crefname{cor}{Corollary}{Corollaries}
\crefname{claim}{Claim}{Claims}
\Crefname{claim}{Claim}{Claims}
\crefname{fact}{Fact}{Facts}
\Crefname{fact}{Fact}{Facts}
\crefname{conj}{Conjecture}{Conjectures}
\Crefname{conj}{Conjecture}{Conjectures}
\crefname{defn}{Definition}{Definitions}
\Crefname{defn}{Definition}{Definitions}
\crefname{exam}{Example}{Examples}
\Crefname{exam}{Example}{Examples}
\crefname{setup}{Setup}{Setups}
\Crefname{setup}{Setup}{Setups}
\crefname{exer}{Exercise}{Exercises}
\Crefname{exer}{Exercise}{Exercises}
\crefname{ques}{Question}{Questions}
\Crefname{ques}{Question}{Questions}
\crefname{rmk}{Remark}{Remarks}
\Crefname{rmk}{Remark}{Remarks}

\counterwithin{figure}{section}

\usepackage{todonotes}

\title{Katona's Shadow Equality Persists under Bounded VC-Dimension}
\author{
Xiaochen Zhao\thanks{School of Mathematical Sciences, Capital Normal University, Beijing, China. Email: 2250501013@cnu.edu.cn.}
\and 
Gennian Ge\thanks{School of Mathematical Sciences, Capital Normal University, Beijing, China. Email: gnge@zju.edu.cn. Gennian Ge is supported by the National Key Research and Development Program of China under Grant 2025YFC3409900, the National Natural Science Foundation of China under Grant 12231014, and Beijing Scholars Program.}
}
\date{}

\begin{document}
\maketitle

\begin{abstract}
Katona's intersecting-shadow theorem states that every
$t$-intersecting $r$-uniform family has size at most
that of its $(r-t)$-shadow, with equality precisely
for complete $r$-uniform families on $2r-t$ points.
We prove that both the shadow inequality and its
equality classification remain valid under the
substantially weaker assumption that the VC-dimension
is at most $r-t$.
By contrast, bounded VC-dimension and $t$-intersection do not, in general, yield the same extremal bounds in the corresponding cardinality problem.
\end{abstract}



\section{Introduction}

The Vapnik--Chervonenkis (VC) dimension is a fundamental concept in statistical learning theory, where it plays a central role in uniform convergence and generalization. It is also an important parameter in extremal combinatorics. For a family $\mathcal F\subseteq 2^{[n]}$ and a set $X\subseteq[n]$, write $\mathcal F|_X=\{F\cap X:F\in\mathcal F\}$. A set $X$ is \emph{shattered} by $\mathcal F$ if $\mathcal F|_X=2^X$, and the VC-dimension \(\operatorname{VCdim}(\mathcal F)\) is the largest size of a shattered set. The classical Sauer--Shelah--Vapnik--Chervonenkis lemma asserts that every family with VC-dimension at most $d$ satisfies $|\mathcal F|\le\sum_{i=0}^d\binom ni$, and this bound is sharp \cite{Sauer1972,Shelah1972,VapnikChervonenkis1971}.

For an $r$-uniform family $\mathcal F$, the classical $t$-intersection condition $|F\cap G|\ge t$ for all $F,G\in\mathcal F$ implies $\operatorname{VCdim}(\mathcal F)\le r-t$. More precisely, putting $k=r-t$, the $t$-intersection condition is equivalent to the absence of the empty trace on every $(k+1)$-set contained in a member of $\mathcal F$. Indeed, two members \(F,G\) violate the \(t\)-intersection condition if and only if \(F\setminus G\) contains a \((k+1)\)-set. By contrast, $\operatorname{VCdim}(\mathcal F)\le k$ requires only that every $(k+1)$-set has at least one missing trace. Hence bounded VC-dimension imposes a substantially weaker restriction than $t$-intersection.

Frankl and Pach initiated the study of the maximum size of a $(d+1)$-uniform family with VC-dimension at most $d$, establishing the upper bound $\binom nd$ \cite{FranklPach1984}. Erd\H{o}s and, independently, Frankl and Pach conjectured that, for sufficiently large $n$, this bound could be sharpened to the Erd\H{o}s--Ko--Rado star bound $\binom{n-1}{d}$ \cite{Erdos1984,EKR1961,FranklPach1984}. In 1997, Ahlswede and Khachatrian disproved the conjecture by constructing families of size $\binom{n-1}{d}+\binom{n-4}{d-2}$ \cite{AhlswedeKhachatrian1997}.
Subsequently, Mubayi and Zhao improved the Frankl--Pach bound for prime-power $d$ and conjectured that the Ahlswede--Khachatrian construction was optimal for sufficiently large $n$ \cite{MubayiZhao2007}. 
More recently, the problem has seen a series of advances concerning both upper bounds and extremal constructions. In particular, the Frankl--Pach bound was shown to be unattainable for every $d\ge2$ and $n\ge2d+2$ \cite{GeXuYipZhangZhao2026}. An asymptotic upper bound matching the leading term of the conjectured construction was obtained in \cite{ChaoXuYipZhang2025}, and was subsequently sharpened to match the order of the second term in \cite{YangYu2025}. The Mubayi--Zhao conjecture was also confirmed for $d=2$ \cite{WangXuZhang2025}. Nevertheless, two recent independent constructions surpassed the Ahlswede--Khachatrian bound for every $d\ge3$, overturning the conjecture of Mubayi and Zhao \cite{TranXu2026,ZhaoGe2026}. These developments show that bounded VC-dimension does not, in general, lead to the same cardinality extremal behavior as the classical intersection condition.

For shadow inequalities, however, a different phenomenon emerges. Given a family $\mathcal F\subseteq2^{[n]}$ and $0\le j\le n$, its $j$-shadow is defined by $\partial_j\mathcal F=\{T\in\binom{[n]}j:T\subseteq F\text{ for some }F\in\mathcal F\}$. A classical result in this direction is Katona's intersecting-shadow theorem.

\begin{theorem}[Katona~\cite{Katona1964}]\label{thm:katona}
Let $1\le t\le r\le n$, and let $\mathcal F\subseteq\binom{[n]}r$ be $t$-intersecting. Then $|\mathcal F|\le|\partial_{r-t}\mathcal F|$.
If $\mathcal F\ne\varnothing$, equality holds if and only if $\mathcal F=\binom Qr$ for some $Q\subseteq[n]$ with $|Q|=2r-t$.
\end{theorem}

Katona's shadow inequality remains valid under bounded VC-dimension: if $\mathcal F\subseteq\binom{[n]}r$ satisfies $\operatorname{VCdim}(\mathcal F)\le k<r$, then $|\mathcal F|\le|\partial_k\mathcal F|$. This follows from the inclusion-matrix independence underlying the
Frankl--Pach method \cite{FranklPach1984}; see also the shadow inequalities of Frankl and F\"uredi \cite{FranklFuredi1984}. Although the inequality was already known, its equality cases had not been characterized under the VC-dimension hypothesis. The remaining problem was to determine whether equality still forces $\mathcal F$ to be a complete $r$-uniform family on $r+k$ points. We show that Katona's equality classification persists under bounded VC-dimension.

\begin{theorem}\label{thm:main}
Let $0\le k<r\le n$, and let $\mathcal F\subseteq\binom{[n]}r$ satisfy $\operatorname{VCdim}(\mathcal F)\le k$. Then $|\mathcal F|\le|\partial_k\mathcal F|$. If $\mathcal F\ne\varnothing$, equality holds if and only if
$\mathcal F=\binom Qr$ for some $Q\subseteq[n]$ with $|Q|=r+k$.
\end{theorem}

The proof hinges on a matrix identity derived from the unique expansions of $(k+1)$-shadow monomials in the degree-$k$ shadow basis
available under equality. This uniqueness also forces each $(k+1)$-set in the shadow to have exactly one missing trace. The identity yields an exact missing-trace descent, while a polynomial argument rules out singleton missing traces. Consequently, all missing traces are empty, forcing $\mathcal F$ to be $(r-k)$-intersecting. The equality classification then follows from Katona's theorem.

Theorem~\ref{thm:main} also determines exactly when the ambient Frankl--Pach bound can be attained.

\begin{cor}\label{cor:ambient}
Let $1\le k<r\le n$, and let $\mathcal F\subseteq\binom{[n]}r$ satisfy $\operatorname{VCdim}(\mathcal F)\le k$. Then $|\mathcal F|\le\binom nk$, with equality if and only if $n=r+k$ and $\mathcal F=\binom{[n]}r$. In particular, $|\mathcal F|<\binom nk$ whenever $n\ne r+k$.
\end{cor}

The corollary follows immediately from $|\mathcal F|\le|\partial_k\mathcal F|\le\binom nk$ and Theorem~\ref{thm:main}. Taking $r=k+1$ recovers the strict Frankl--Pach bound for $n\ge2k+2$ \cite{GeXuYipZhangZhao2026}, while the general statement determines the equality cases for every $1\le k<r$.

\noindent\textbf{Organization of the paper.}
Section~2 develops the algebraic framework for shadow equality and derives the key matrix identity, while Section~3 derives $(r-k)$-intersection and completes the classification via Katona's theorem. Section~4 proposes an extension to Sperner families, relates it to Frankl's conjecture, and establishes the case for $d\le2$, whose proof is given in Appendix~A.

\section{The equality normal form}\label{sec:equality-normal-form}
For $A\subseteq[n]$, write $x_A=\prod_{i\in A}x_i$, where $x_\varnothing=1$. We identify each $F\subseteq[n]$ with its characteristic vector $v_F\in\{0,1\}^n$. A polynomial restricted to $\mathcal F$ will be regarded as a function on the points $\{v_F:F\in\mathcal F\}$; in particular, all polynomial identities below are understood on this restricted domain.

\begin{prop}\label{prop:basis}
Let $0\le k<r\le n$, and suppose $\mathcal F\subseteq\binom{[n]}r$ satisfies $\operatorname{VCdim}(\mathcal F)\le k$. Then $|\mathcal F|\le|\partial_k\mathcal F|$. If equality holds, the inclusion matrix $W_k=(\mathbf1_{\{T\subseteq F\}})_{F\in\mathcal F,\,T\in\partial_k\mathcal F}$ is square and nonsingular. Equivalently, the functions $\{x_T|_{\mathcal F}:T\in\partial_k\mathcal F\}$ form a basis of $\mathbb R^{\mathcal F}$.
\end{prop}

\begin{proof}
We first show that the monomials of degree at most $k$ span $\mathbb R^{\mathcal F}$. For each $(k+1)$-set $Y\subseteq[n]$, the VC-dimension hypothesis supplies a missing trace $B_Y\subseteq Y$. Consequently, $q_Y=x_{B_Y}\prod_{a\in Y\setminus B_Y}(x_a-1)$ vanishes on $\mathcal F$. Its unique term of degree $k+1$ is $x_Y$, with coefficient one, so this relation expresses $x_Y$ on $\mathcal F$ as a linear combination of monomials of degree at most $k$. More generally, if $|W|>k+1$, choose $Y\in\binom W{k+1}$, multiply the relation for $x_Y$ by $x_{W\setminus Y}$, and multilinearize. The leading term becomes $x_W$, while every remaining term has degree smaller than $|W|$. Induction on $|W|$ reduces all higher-degree monomials. Since squarefree monomials span the functions on $\{0,1\}^n$, their restrictions of degree at most $k$ span $\mathbb R^{\mathcal F}$.

The uniformity assumption allows us to pass from degrees \emph{at most} $k$ to degree \emph{exactly} $k$. For $|A|=j\le k$ we have, on the $r$-slice,
\[
 x_A=\binom{r-j}{k-j}^{-1}\sum_{\substack{T\in\binom{[n]}k\\A\subseteq T}}x_T.
\]
Indeed, for an $r$-set $F$, both sides vanish unless $A\subseteq F$; in that case, exactly $\binom{r-j}{k-j}$
of the indexed $k$-sets are contained in $F$. It follows that the degree-$k$ monomials span $\mathbb R^{\mathcal F}$. Those indexed by
$T\notin\partial_k\mathcal F$ vanish on $\mathcal F$, so the remaining $|\partial_k\mathcal F|$ functions still span. This proves the inequality, and equality makes this spanning family a basis.
\end{proof}

For the rest of the equality analysis, assume that $\mathcal F\ne\varnothing$, $1\le k<r$, and $|\mathcal F|=|\partial_k\mathcal F|$; the case $k=0$ will be handled separately in the proof of the main theorem. Write $S_k=\partial_k\mathcal F$ and $S_{k+1}=\partial_{k+1}\mathcal F$. For $E\in S_{k+1}$ and a missing trace $B\subsetneq E$, consider $q_{E,B}=x_B\prod_{a\in E\setminus B}(x_a-1)$ and $\psi_{E,B}=x_E-q_{E,B}$. The first polynomial vanishes on $\mathcal F$, and the degree-$(k+1)$ terms cancel in the second. Thus $\psi_{E,B}$ has degree at most $k$ and represents $x_E|_{\mathcal F}$. We next compute its coordinates in the shadow basis from \cref{prop:basis}.

\begin{lemma}\label{lem:coeff}
Let $m=|B|$, $s=k-m$, and denote by $c_B(T,E)$ the coefficient of $x_T|_{\mathcal F}$ in the shadow-basis expansion of $\psi_{E,B}|_{\mathcal F}$. If $B\not\subseteq T$, then $c_B(T,E)=0$. If $B\subseteq T$ and $t=|T\cap(E\setminus B)|$, then
\begin{equation}\label{eq:coeff}
 c_B(T,E)=
 \frac{(-1)^{s-t}\binom{r-k+t-1}{t}}
 {\binom{r-m}{k-m}\binom{k-m}{t}}.
\end{equation}
In particular, for $B\subseteq T\subsetneq E$ with $|T|=k$, we have $c_B(T,E)=(r-k)/(r-m)>0$.
\end{lemma}

\begin{proof}
Set $L=E\setminus B$, so $|L|=k+1-m=s+1$. Expanding $q_{E,B}$ and canceling $x_E$ gives
\[
 \psi_{E,B}=-\sum_{A\subsetneq L}(-1)^{|L|-|A|}x_{B\cup A}.
\]
Every monomial contains $B$, and the homogenization identity in \cref{prop:basis} only replaces a monomial by those indexed by supersets. Hence $c_B(T,E)=0$ when $B\not\subseteq T$. Suppose $B\subseteq T$. For a term $x_{B\cup A}$ to contribute to the coefficient at $x_T$, we must have $A\subseteq T\cap L$. Writing $t=|T\cap L|$, homogenization therefore gives
\[
 c_B(T,E)=\sum_{i=0}^{t}(-1)^{s-i}\binom ti
 \binom{r-m-i}{k-m-i}^{-1}.
\]
The sign $(-1)^{s-i}$ follows from $-(-1)^{|L|-i}=(-1)^{s-i}$.

To evaluate the sum, write $R=r-m$. The binomial identity
\[
 \frac{\binom ti}{\binom{R-i}{s-i}}
 =\frac{\binom Ri\binom{s-i}{t-i}}
 {\binom Rs\binom st}
\]
reduces the calculation to
\[
 c_B(T,E)=\frac{(-1)^s}{\binom Rs\binom st}
 \sum_{i=0}^t(-1)^i\binom Ri\binom{s-i}{t-i}.
\]
The remaining alternating sum can be evaluated by coefficient extraction. Indeed,
\[
 \begin{aligned}
 \sum_{i=0}^t(-1)^i\binom Ri\binom{s-i}{t-i}
 &=[u^t]\sum_{i=0}^R(-1)^i\binom Ri
       u^i(1+u)^{s-i}\\
 &=[u^t](1+u)^s\left(1-\frac{u}{1+u}\right)^R\\
 &=[u^t](1+u)^{s-R}
 =(-1)^t\binom{R-s+t-1}{t}.
 \end{aligned}
\]
Substituting $R=r-m$ and $R-s=r-k$ proves \eqref{eq:coeff}. If $B\subseteq T\subsetneq E$, then $t=k-m=s$; the formula simplifies to $(r-k)/(r-m)$.
\end{proof}

\begin{claim}
Every $E\in S_{k+1}$ has exactly one missing trace $b(E)\subsetneq E$.
\end{claim}
\begin{proof}
The VC-dimension bound ensures that every $E\in S_{k+1}$ has a missing trace, while its full trace occurs because some $F\in\mathcal F$ contains $E$. Suppose $B,B'\subsetneq E$ are both missing. The polynomials $\psi_{E,B}$ and $\psi_{E,B'}$ represent the same function $x_E|_{\mathcal F}$. Since the shadow monomials form a basis, their coordinates agree. For each $e\in E$, the facet $T_e=E\setminus\{e\}$ belongs to $S_k$, and \cref{lem:coeff} gives
\[
 c_B(T_e,E)=
 \begin{cases}
 0,&e\in B,\\
 \dfrac{r-k}{r-|B|},&e\notin B.
 \end{cases}
\]
Thus the zero coordinates among the facets of $E$ record precisely the elements of $B$. The equality of the two coordinate vectors forces $B=B'$, proving uniqueness.
\end{proof}

Define
\begin{equation*}
  C=(c(T,E))_{T\in S_k,\,E\in S_{k+1}}\quad\text{by}\quad c(T,E)=c_{b(E)}(T,E),
\end{equation*}
and let 
\[
 U=(\mathbf1_{\{T\subseteq E\}})_{E\in S_{k+1},\,T\in S_k}.
\]

\begin{prop}\label{prop:normal}
\begin{equation}\label{eq:normal}
 CU=(r-k)I_{S_k}.
\end{equation}
\end{prop}

\begin{proof}
The matrix identity comes from expressing the same incidence matrix in two ways. Let $W_{k+1}=(\mathbf1_{\{E\subseteq F\}})_{F\in\mathcal F,\,E\in S_{k+1}}$. The expansions furnished by the unique missing traces give $W_{k+1}=W_kC$. On the other hand, for fixed $F\in\mathcal F$ and $T\in S_k$, the $(F,T)$-entry of $W_{k+1}U$ counts $(k+1)$-sets $E$ with $T\subseteq E\subseteq F$. This number equals $r-k$ when $T\subseteq F$, and zero otherwise. Hence $W_{k+1}U=(r-k)W_k$. Substituting the first identity and canceling the nonsingular matrix $W_k$ yields \eqref{eq:normal}.
\end{proof}

The matrix identity \eqref{eq:normal} imposes local constraints on missing traces that will be crucial to the descent argument in Section~3. Specifically, its diagonal entries yield a weighted identity that guarantees the existence of certain extensions, while its off-diagonal entries impose a compatibility condition on the missing traces of adjacent extensions.

\begin{lemma}\label{lem:matrix-entries}
For $T\in S_k$ and $z\notin T$, put $E_z=T\cup\{z\}$, and define $\Gamma(T)=\{z\notin T:E_z\in S_{k+1},\ b(E_z)\subseteq T\}$.
Then the following hold.
\begin{enumerate}[(i)]
\item We have $\sum_{z\in\Gamma(T)}(r-|b(E_z)|)^{-1}=1$. In particular, $\Gamma(T)\ne\varnothing$.
\item For $h\in T$ and $y\notin T$ with $E_y\in S_{k+1}$, we have $h\in b(E_y)$ if and only if $h\in b(E_z)$ for every $z\in\Gamma(T)\setminus\{y\}$.
\end{enumerate}
\end{lemma}

\begin{proof}
For $T_1,T_2\in S_k$, the $(T_1,T_2)$-entry of \eqref{eq:normal} is the sum of $c(T_1,E)$ over all $E\in S_{k+1}$ containing $T_2$, and equals
$(r-k)\mathbf1_{\{T_1=T_2\}}$. By \cref{lem:coeff}, $c(T_1,E)$ vanishes unless $b(E)\subseteq T_1$, and every nonzero coefficient has sign
$(-1)^{k-|T_1\cap E|}$.

For the diagonal entry indexed by $T$, the nonzero summands correspond to these $E\in S_{k+1}$ with $b(E)\subseteq T \subseteq E$, precisely to $z\in\Gamma(T)$, each equal to $(r-k)/(r-|b(E_z)|)$ by \cref{lem:coeff}. Dividing the diagonal identity by $r-k$ proves (i).

For (ii), put $T'=(T\setminus\{h\})\cup\{y\}\in S_k$ and consider the $(T',T)$-entry of \eqref{eq:normal}, which is zero. Since $T'=E_y\setminus\{h\}$, \cref{lem:coeff} gives $c(T',E_y)\ge0$, with equality precisely when $h\in b(E_y)$. For $z\ne y$, we have
$T'\cap E_z=T\setminus\{h\}$, so $c(T',E_z)\le0$, with strict inequality precisely when $b(E_z)\subseteq T'$, which is equivalent to $z\in\Gamma(T)$ and $h\notin b(E_z)$. Since the sum of these coefficients is zero, the positive term vanishes if and only if all negative terms vanish. This proves (ii).
\end{proof}

\section{Descent and rigidity}\label{sec:descent}
We first establish an exact descent for missing traces using \cref{lem:matrix-entries}. We then derive a rigidity statement for singleton missing traces and use it to rule them out.

\begin{prop}\label{prop:descent}
Let $E\in S_{k+1}$ have missing trace $B=b(E)$ of size $m\ge1$. For every $x\in B$, there are exactly $r-m+1$
points $u\notin E$ such that $(E\setminus\{x\})\cup\{u\}\in S_{k+1}$ and $b\bigl((E\setminus\{x\})\cup\{u\}\bigr)=B\setminus\{x\}$.
\end{prop}

\begin{proof}
Fix $x\in B$ and put $R=E\setminus\{x\}$. Let $\Lambda=\Gamma(R)$, where $\Gamma$ is defined in \cref{lem:matrix-entries}. For each $u\in\Lambda$, write $E^u=R\cup\{u\}$ and $B_u=b(E^u)\subseteq R$. Since $x\in b(E)$, we have $x\notin\Lambda$, so every point of $\Lambda$ lies outside $E$. By \cref{lem:matrix-entries}(i), $\Lambda\ne\varnothing$.

We first show that all the sets $B_u$, $u\in\Lambda$, coincide. Fix $u\in\Lambda$ and $a\in B_u=b(E^u)$. Applying \cref{lem:matrix-entries}(ii) with $T=R$, $h=a$, and $E_y=E^u$, we obtain $a\in B_v$ for every $v\in\Lambda\setminus\{u\}$. Since $u$ and $a$ were arbitrary, all the sets $B_u$, $u\in\Lambda$, coincide; let $B^*$ denote their common value.

To identify $B^*$, fix $a\in R$. Applying \cref{lem:matrix-entries}(ii) with $T=R$, $h=a$, and $E_y=E$, and using $x\notin\Lambda$, we obtain
\[
a\in B \quad\Longleftrightarrow\quad a\in B_u\text{ for every }u\in\Lambda.
\]
Since $\Lambda\ne\varnothing$ and all the $B_u$ equal $B^*$, this is equivalent to $a\in B$ if and only if $a\in B^*$. As $B^*\subseteq R$, we conclude that $B^*=B\setminus\{x\}$.

Finally, \cref{lem:matrix-entries}(i) gives
\[
\frac{|\Lambda|}{r-|B^*|}=1,
\]
so $|\Lambda|=r-m+1$. Thus, exactly \(r-m+1\) points outside \(E\) have the required properties.
\end{proof}

Iterating the above descent proposition, we obtain a singleton missing trace whenever a nonempty one exists. The following lemma establishes strong rigidity properties forced by the existence of such a singleton, which will allow us to rule out this possibility.

\begin{lemma}\label{lem:fan}
Suppose $L\in\binom{[n]}k$ and $y\notin L$ satisfy $L\cup\{y\}\in S_{k+1}$ and $b(L\cup\{y\})=\{y\}$. Then there exists $Z\in\mathcal F$ with the following properties:
\begin{enumerate}[(i)]
\item $Z\cap(L\cup\{y\})=\varnothing$, and $Z$ is the unique member of $\mathcal F$ disjoint from $L$;
\item $L\cup\{z\}\in S_{k+1}$ and $b(L\cup\{z\})=\varnothing$ for every $z\in Z$;
\item writing $Q_0=L\cup Z$, if $p\notin Q_0$ and $E\in S_{k+1}$ contains $p$, then $b(E)\not\subseteq L$.
\end{enumerate}
\end{lemma}

\begin{proof}
Proposition~\ref{prop:descent}, applied to the singleton missing trace $\{y\}$, gives exactly $r$ points $z_1,\ldots,z_r\notin L\cup\{y\}$ for which $b(L\cup\{z_i\})=\varnothing$. Because $\{y\}$ is the \emph{only} missing trace on $L\cup\{y\}$, the empty trace occurs. Choose $Z\in\mathcal F$ disjoint from $L\cup\{y\}$. For each $i$, the missing empty trace on $L\cup\{z_i\}$ means that every member of $\mathcal F$ intersects $L\cup\{z_i\}$. Since $Z\cap L=\varnothing$, it must contain $z_i$. Thus $\{z_1,\ldots,z_r\}\subseteq Z$, and $|Z|=r$ forces equality. Any other member disjoint from $L$ would likewise contain all $z_i$ and hence coincide with $Z$. This establishes (i) and (ii).

The proof of (iii) requires an additional use of the shadow basis. Consider $P_L=\prod_{\ell\in L}(1-x_\ell)$. On the Boolean cube this polynomial is the indicator of sets disjoint from $L$. Since $Z$ is the unique member of $\mathcal F$ disjoint from $L$, its restriction to
$\mathcal F$ is the point indicator $\delta_Z$. Write the unique shadow-basis expansion as
\[
 P_L|_{\mathcal F}=\sum_{T\in S_k}a_Tx_T|_{\mathcal F}.
\]
We need the signs of $a_T$. Put $q=|T\cap L|$. Expanding $P_L$ and then homogenizing each $x_A$ gives
\[
 a_T=\sum_{i=0}^{q}(-1)^i\binom qi \binom{r-i}{k-i}^{-1}
 =(-1)^q\frac{\binom{r-k+q-1}{q}}{\binom rk\binom kq}.
\]
The last identity follows by the same alternating-sum evaluation used in \cref{lem:coeff}. In particular, $a_T$ is nonzero and has sign $(-1)^{|T\cap L|}$.

Now fix $p\notin Q_0=L\cup Z$. Since $P_L|_{\mathcal F}=\delta_Z$ and $p\notin Z$, we have $x_pP_L=0$ on $\mathcal F$. Multiplying the shadow-basis expansion of $P_L$ by $x_p$ and using the expansions from \cref{lem:coeff}, we obtain
\[
\begin{aligned}
0=x_pP_L|_{\mathcal F}
&=\sum_{\substack{T\in S_k\\p\in T}} a_Tx_T|_{\mathcal F}
\\&\quad+
\sum_{\substack{T\in S_k\\p\notin T\\T\cup\{p\}\in S_{k+1}}}a_T\sum_{U\in S_k}c(U,T\cup\{p\})x_U|_{\mathcal F}.
\end{aligned}
\]
Here we have omitted the terms with $T\cup\{p\}\notin S_{k+1}$, which vanish on $\mathcal F$.

Since $p\notin L$, the first sum has no contribution to the $x_L$-coordinate. Comparing this coordinate on both sides gives
\[
0=\sum_{\substack{T\in S_k\\p\notin T\\T\cup\{p\}\in S_{k+1}}}a_Tc(L,T\cup\{p\}).
\]
For every nonzero summand, the coefficient formulas above and in \cref{lem:coeff} give
\[
\operatorname{sgn}(a_T)=(-1)^{|T\cap L|},\quad \operatorname{sgn}(c(L,T\cup\{p\}))=(-1)^{k-|T\cap L|},
\]
where we used $p\notin L$. Thus every nonzero summand has sign $(-1)^k$. Since their sum is zero, all summands must vanish. As $a_T\ne0$ for every $T\in S_k$, we obtain $c(L,T\cup\{p\})=0$ for every indexed $T$. By \cref{lem:coeff}, this means that $b(E)\not\subseteq L$ for every
$E\in S_{k+1}$ containing $p$, proving (iii).
\end{proof}

\begin{claim}\label{claim:empty}
For every $E\in S_{k+1}$, we have $b(E)=\varnothing$.
\end{claim}

\begin{proof}
Suppose that some $E\in S_{k+1}$ has a nonempty missing trace. By repeated applications of \cref{prop:descent}, we obtain $E_1=L\cup\{y\}\in S_{k+1}$ with $b(E_1)=\{y\}$. Applying \cref{lem:fan} to $E_1$, we obtain $Z\in\mathcal F$ disjoint from $E_1$. Put $Q=L\cup Z$, so $|Q|=r+k$.
Moreover, \cref{lem:fan}(iii) implies that
\begin{equation}\label{eq:empty-face-Q}
  D\in S_{k+1},\quad b(D)=\varnothing\quad\Longrightarrow\quad D\subseteq Q.
\end{equation}

Let $H\in\mathcal F$ be some $r$-set containing $E_1$. Since $Z\cap E_1=\varnothing$ and $|H|=|Z|=r$, we have $|Z\setminus H|=|H\setminus Z|\ge k+1$. Choose $D\in\binom{Z\setminus H}{k+1}$. Then $Z$ and $H$ realize the full and empty traces on $D$, respectively. Since $D\in S_{k+1}$ has a unique missing trace, we have $\varnothing\ne b(D)\subsetneq D$. Fix $u\in b(D)$. Applying \cref{prop:descent} repeatedly to
the elements of $b(D)\setminus\{u\}$, we obtain $E_2=M\cup\{u\}\in S_{k+1}$ with $b(E_2)=\{u\}$.

Now apply \cref{lem:fan} to $E_2$. This gives $W\in\mathcal F$ disjoint from $E_2$ such that $b(M\cup\{w\})=\varnothing$ for every $w\in W$.
By \eqref{eq:empty-face-Q}, we have $M\cup\{w\}\subseteq Q$ for every $w\in W$, and hence $M\cup W\subseteq Q$. On the other hand, $u\notin M\cup W$, because $u\notin M$ and $W\cap E_2=\varnothing$. Therefore
\[
  M\cup W\subseteq Q\setminus\{u\}.
\]
This is impossible, since $u\in D\subseteq Z\subseteq Q$, $M\cap W=\varnothing$, and $|M\cup W|=k+r=|Q|$. Thus every $E\in S_{k+1}$ has
$b(E)=\varnothing$.
\end{proof}

\begin{proof}[Proof of \cref{thm:main}]
The shadow inequality follows from \cref{prop:basis}. We now characterize the equality cases.

Suppose $\mathcal F\ne\varnothing$ and $|\mathcal F|=|\partial_k\mathcal F|$. If $k=0$, then $\partial_0\mathcal F=\{\varnothing\}$, so $\mathcal F$ consists of a single $r$-set, as required.

Assume $k\ge1$. By \cref{claim:empty}, the empty trace is missing from every $E\in S_{k+1}$. If some $F,G\in\mathcal F$ satisfied $|F\setminus G|\ge k+1$, then any $E\in\binom{F\setminus G}{k+1}$ would belong to $S_{k+1}$, while $G\cap E=\varnothing$ would realize
its empty trace, a contradiction. Thus $|F\cap G|\ge r-k$ for all $F,G\in\mathcal F$. Hence $\mathcal F$ is $(r-k)$-intersecting,
and the equality case of \cref{thm:katona} gives $\mathcal F=\binom Qr$ for some $Q\subseteq[n]$ with $|Q|=r+k$.

Conversely, let $\mathcal F=\binom Qr$ for some $Q\subseteq[n]$ with $|Q|=r+k$. Then $\mathcal F$ is $(r-k)$-intersecting, and hence $\operatorname{VCdim}(\mathcal F)\le k$. Moreover, $\partial_k\mathcal F=\binom Qk$, so
\[
|\mathcal F|=\binom{r+k}{r}=\binom{r+k}{k}=|\partial_k\mathcal F|.
\]
This completes the proof.
\end{proof}

\section{Concluding remarks}\label{sec:concluding}

An appealing extension of Theorem~\ref{thm:main} is to replace uniformity by the Sperner property. This retains the incomparability of distinct members while allowing their sizes to vary, raising the question of whether the shadow inequality and its equality classification persist without the fixed-size structure underlying our proof.

\begin{conj}\label{conj:antichain-shadow}
Let $d\ge1$, and let $\mathcal A\subseteq2^{[n]}$ be a nonempty antichain with $\operatorname{VCdim}(\mathcal A)\le d$ and $|A|\ge d+1$ for every $A\in\mathcal A$. Then $|\mathcal A|\le|\partial_d\mathcal A|$. Moreover, equality holds if and only if $\mathcal A=\binom Qr$ for some $Q\subseteq[n]$ with $r>d$ and $|Q|=r+d$.
\end{conj}

This would strengthen a classical conjecture of Frankl concerning Sperner families of bounded VC-dimension.

\begin{conj}[Frankl~\cite{Frankl1989}]\label{conj:frankl}
Let $n\ge2d\ge2$, and let $\mathcal F\subseteq2^{[n]}$ be an antichain with $\operatorname{VCdim}(\mathcal F)\le d$. Then $|\mathcal F|\le\binom nd$.
\end{conj}

Frankl proved this conjecture for $d\le2$, and Anstee and Sali established the case $d=3$ \cite{Frankl1989,AnsteeSali1997}. More recently, Ge, Wang, Xu, and Zhang proved the conjecture for every fixed $d$ and all sufficiently large $n$, and characterized the equality cases \cite{GeWangXuZhang2026}. The full range $n\ge2d$ remains open. The connection with \cref{conj:antichain-shadow} is made precise by the following implication.

\begin{prop}\label{prop:shadow-implies-frankl}
\cref{conj:antichain-shadow} implies \cref{conj:frankl}.
\end{prop}

\begin{proof}
Write $\mathcal F=\mathcal S\,\cup\,\mathcal L$, where $\mathcal S=\{F\in\mathcal F:|F|\le d\}$ and $\mathcal L=\{F\in\mathcal F:|F|\ge d+1\}$. Since $n\ge2d$ and \(\mathcal S\) is an antichain, a symmetric-chain decomposition of $2^{[n]}$ gives an injection $\phi:\mathcal S\to\binom{[n]}d$ satisfying $S\subseteq\phi(S)$ for every $S\in\mathcal S$. Moreover, $\phi(\mathcal S)\cap\partial_d\mathcal L=\varnothing$,
since otherwise $S\subseteq\phi(S)\subseteq L$ for some $S\in\mathcal S$ and $L\in\mathcal L$, contradicting the Sperner property. Hence $|\mathcal S|\le\binom nd-|\partial_d\mathcal L|$. If $\mathcal L=\varnothing$, the desired bound follows immediately. Otherwise, applying \cref{conj:antichain-shadow} to $\mathcal L$ gives $|\mathcal L|\le|\partial_d\mathcal L|$. Combining these inequalities yields $|\mathcal F|\le\binom nd$.
\end{proof}

The stronger shadow statement can be established for $d\le2$, including its complete equality classification.

\begin{theorem}\label{thm:antichain-d2}
Let $d\in\{1,2\}$, and let $\mathcal A\subseteq2^{[n]}$ be a nonempty antichain with $\operatorname{VCdim}(\mathcal A)\le d$ and $|A|\ge d+1$ for every $A\in\mathcal A$. Then $|\mathcal A|\le|\partial_d\mathcal A|$, with equality if and only if $\mathcal A=\binom Qr$ for some $Q\subseteq[n]$ with $r>d$ and $|Q|=r+d$.
\end{theorem}

A self-contained proof is given in Appendix~\ref{app:antichain-d2}. Resolving \cref{conj:antichain-shadow} in full may require new ideas for handling the interaction between different member sizes.

\bibliographystyle{abbrv}
\bibliography{references}

@article {AhlswedeKhachatrian1997,
    AUTHOR = {Ahlswede, R. and Khachatrian, L. H.},
     TITLE = {Counterexample to the {F}rankl-{P}ach conjecture for uniform,
              dense families},
   JOURNAL = {Combinatorica},
  FJOURNAL = {Combinatorica. An International Journal on Combinatorics and
              the Theory of Computing},
    VOLUME = {17},
      YEAR = {1997},
    NUMBER = {2},
     PAGES = {299--301},
      ISSN = {0209-9683,1439-6912},
   MRCLASS = {05D05},
  MRNUMBER = {1479304},
MRREVIEWER = {Martin\ Klazar},
       DOI = {10.1007/BF01200912},
       URL = {https://doi.org/10.1007/BF01200912},
}

@article {AnsteeSali1997,
    AUTHOR = {Anstee, R. P. and Sali, A.},
     TITLE = {Sperner families of bounded {VC}-dimension},
   JOURNAL = {Discrete Math.},
  FJOURNAL = {Discrete Mathematics},
    VOLUME = {175},
      YEAR = {1997},
    NUMBER = {1-3},
     PAGES = {13--21},
      ISSN = {0012-365X,1872-681X},
   MRCLASS = {05D05 (05C65 06A07)},
  MRNUMBER = {1475834},
MRREVIEWER = {Konrad\ Engel},
       DOI = {10.1016/S0012-365X(96)00137-9},
       URL = {https://doi.org/10.1016/S0012-365X(96)00137-9},
}

@article {ChaoXuYipZhang2025,
    AUTHOR = {Chao, Ting-Wei and Xu, Zixiang and Yip, Chi Hoi and Zhang,
              Shengtong},
     TITLE = {Uniform set systems with small {VC}-dimension},
   JOURNAL = {Int. Math. Res. Not. IMRN},
  FJOURNAL = {International Mathematics Research Notices. IMRN},
      YEAR = {2025},
    NUMBER = {17},
     PAGES = {Paper No. rnaf269, 20},
      ISSN = {1073-7928,1687-0247},
   MRCLASS = {05D05},
  MRNUMBER = {4954468},
MRREVIEWER = {Norihide\ Tokushige},
       DOI = {10.1093/imrn/rnaf269},
       URL = {https://doi.org/10.1093/imrn/rnaf269},
}

@incollection {Erdos1984,
    AUTHOR = {Erd\H{o}s, Paul},
     TITLE = {On some problems in graph theory, combinatorial analysis and
              combinatorial number theory},
 BOOKTITLE = {Graph theory and combinatorics ({C}ambridge, 1983)},
     PAGES = {1--17},
 PUBLISHER = {Academic Press, London},
      YEAR = {1984},
      ISBN = {0-12-111760-X},
   MRCLASS = {05-02 (11Bxx)},
  MRNUMBER = {777160},
MRREVIEWER = {Stefan\ A.\ Burr},
}

@article {EKR1961,
    AUTHOR = {Erd\H{o}s, P. and Ko, Chao and Rado, R.},
     TITLE = {Intersection theorems for systems of finite sets},
   JOURNAL = {Quart. J. Math. Oxford Ser. (2)},
  FJOURNAL = {The Quarterly Journal of Mathematics. Oxford. Second Series},
    VOLUME = {12},
      YEAR = {1961},
     PAGES = {313--320},
      ISSN = {0033-5606,1464-3847},
   MRCLASS = {04.60},
  MRNUMBER = {140419},
MRREVIEWER = {S.\ Ginsburg},
       DOI = {10.1093/qmath/12.1.313},
       URL = {https://doi.org/10.1093/qmath/12.1.313},
}

@article {Frankl1989,
    AUTHOR = {Frankl, P.},
     TITLE = {Traces of antichains},
   JOURNAL = {Graphs Combin.},
  FJOURNAL = {Graphs and Combinatorics},
    VOLUME = {5},
      YEAR = {1989},
    NUMBER = {1},
     PAGES = {295--299},
      ISSN = {0911-0119,1435-5914},
   MRCLASS = {99-04},
  MRNUMBER = {1554355},
       DOI = {10.1007/BF01788682},
       URL = {https://doi.org/10.1007/BF01788682},
}

@incollection {FranklFuredi1984,
    AUTHOR = {Frankl, P. and F\"uredi, Z.},
     TITLE = {Families of finite sets with missing intersections},
 BOOKTITLE = {Finite and infinite sets, {V}ol.\ {I}, {II} ({E}ger, 1981)},
    SERIES = {Colloq. Math. Soc. J\'anos Bolyai},
    VOLUME = {37},
     PAGES = {305--318},
 PUBLISHER = {North-Holland, Amsterdam},
      YEAR = {1984},
      ISBN = {0-444-86763-5},
   MRCLASS = {05A05 (04A20)},
  MRNUMBER = {818242},
MRREVIEWER = {J.\ R.\ Griggs},
}

@article {FranklPach1984,
    AUTHOR = {Frankl, Peter and Pach, J\'anos},
     TITLE = {On disjointly representable sets},
   JOURNAL = {Combinatorica},
  FJOURNAL = {Combinatorica. An International Journal of the J\'anos Bolyai
              Mathematical Society},
    VOLUME = {4},
      YEAR = {1984},
    NUMBER = {1},
     PAGES = {39--45},
      ISSN = {0209-9683},
   MRCLASS = {05A20},
  MRNUMBER = {739411},
MRREVIEWER = {E.\ C.\ Milner},
       DOI = {10.1007/BF02579155},
       URL = {https://doi.org/10.1007/BF02579155},
}

@article {GeXuYipZhangZhao2026,
    AUTHOR = {Ge, Gennian and Xu, Zixiang and Yip, Chi Hoi and Zhang,
              Shengtong and Zhao, Xiaochen},
     TITLE = {The {F}rankl-{P}ach upper bound is not tight for any
              uniformity},
   JOURNAL = {J. Combin. Theory Ser. A},
  FJOURNAL = {Journal of Combinatorial Theory. Series A},
    VOLUME = {217},
      YEAR = {2026},
     PAGES = {Paper No. 106078, 9},
      ISSN = {0097-3165,1096-0899},
   MRCLASS = {05D05},
  MRNUMBER = {4916840},
MRREVIEWER = {William\ Linz},
       DOI = {10.1016/j.jcta.2025.106078},
       URL = {https://doi.org/10.1016/j.jcta.2025.106078},
}

@article {Katona1964,
    AUTHOR = {Katona, Gy.},
     TITLE = {Intersection theorems for systems of finite sets},
   JOURNAL = {Acta Math. Acad. Sci. Hungar.},
  FJOURNAL = {Acta Mathematica. Academiae Scientiarum Hungaricae},
    VOLUME = {15},
      YEAR = {1964},
     PAGES = {329--337},
      ISSN = {0001-5954,1588-2632},
   MRCLASS = {04.60},
  MRNUMBER = {168468},
MRREVIEWER = {H.\ Jerome\ Keisler},
       DOI = {10.1007/BF01897141},
       URL = {https://doi.org/10.1007/BF01897141},
}

@article {MubayiZhao2007,
    AUTHOR = {Mubayi, Dhruv and Zhao, Yi},
     TITLE = {On the {VC}-dimension of uniform hypergraphs},
   JOURNAL = {J. Algebraic Combin.},
  FJOURNAL = {Journal of Algebraic Combinatorics. An International Journal},
    VOLUME = {25},
      YEAR = {2007},
    NUMBER = {1},
     PAGES = {101--110},
      ISSN = {0925-9899,1572-9192},
   MRCLASS = {05C35 (05C65)},
  MRNUMBER = {2284149},
MRREVIEWER = {Jen\H o\ Lehel},
       DOI = {10.1007/s10801-006-0025-4},
       URL = {https://doi.org/10.1007/s10801-006-0025-4},
}

@article {Sauer1972,
    AUTHOR = {Sauer, N.},
     TITLE = {On the density of families of sets},
   JOURNAL = {J. Combinatorial Theory Ser. A},
  FJOURNAL = {Journal of Combinatorial Theory. Series A},
    VOLUME = {13},
      YEAR = {1972},
     PAGES = {145--147},
      ISSN = {0097-3165},
   MRCLASS = {02H10 (04A20 05A05)},
  MRNUMBER = {307902},
MRREVIEWER = {L.\ Bukovsk\'y},
       DOI = {10.1016/0097-3165(72)90019-2},
       URL = {https://doi.org/10.1016/0097-3165(72)90019-2},
}

@article {Shelah1972,
    AUTHOR = {Shelah, Saharon},
     TITLE = {A combinatorial problem; stability and order for models and
              theories in infinitary languages},
   JOURNAL = {Pacific J. Math.},
  FJOURNAL = {Pacific Journal of Mathematics},
    VOLUME = {41},
      YEAR = {1972},
     PAGES = {247--261},
      ISSN = {0030-8730,1945-5844},
   MRCLASS = {02H10 (02B25 04A20 05A05)},
  MRNUMBER = {307903},
MRREVIEWER = {M.\ Makkai},
       URL = {http://projecteuclid.org/euclid.pjm/1102968432},
}

@incollection {VapnikChervonenkis1971,
    AUTHOR = {Vapnik, V. N. and Chervonenkis, A. Ya.},
     TITLE = {On the uniform convergence of relative frequencies of events
              to their probabilities},
 BOOKTITLE = {Measures of complexity},
     PAGES = {11--30},
      NOTE = {Reprint of Theor. Probability Appl. {\bf 16} (1971), 264--280},
 PUBLISHER = {Springer, Cham},
      YEAR = {2015},
      ISBN = {978-3-319-21851-9; 978-3-319-21852-6},
   MRCLASS = {60F05 (60F15)},
  MRNUMBER = {3408730},
}

@article{TranXu2026,
  title={Beating the Ahlswede--Khachatrian bound for the {Erd\H{o}s--Frankl--Pach} problem},
  author={T. Tran and Z. Xu},
  journal={arXiv preprint},
  year={arXiv:2606.23469, 2026},
}

@article{WangXuZhang2025,
  title={Largest 3-uniform set systems with VC-dimension 2},
  author={J. Wang and Z. Xu and S. Zhang},
  journal={arXiv preprint},
  year={arXiv:2505.07756, 2025},
}

@article{YangYu2025,
  title={Maxmum size of a uniform family with bounded VC-dimension},
  author={T. Yang and X. Yu},
  journal={arXiv preprint},
  year={arXiv:2508.14334, 2025},
}

@article{ZhaoGe2026,
  title={Recursive lifting beyond the Ahlswede--Khachatrian construction},
  author={X. Zhao and G. Ge},
  journal={arXiv preprint},
  year={arXiv:2607.04858, 2026},
}

@article{GeWangXuZhang2026,
  title={Largest Sperner families with bounded VC-dimension},
  author={Gennian Ge and Hanyu Wang and Zixiang Xu and Yixuan Zhang},
  journal={arXiv preprint},
  year={arXiv:2611.06860, 2026},
}

\begin{appendices}
\section{Proof of Theorem~\ref{thm:antichain-d2}}\label{app:antichain-d2}
Write $\operatorname{supp}\mathcal A=\bigcup\mathcal A$ and $x_S(A)=\mathbf1_{\{S\subseteq A\}}$. We use the elementary VC interpolation fact that, when $\operatorname{VCdim}(\mathcal A)\le d$, the functions $x_S$ with $|S|\le d$ span $\mathbb R^{\mathcal A}$. Indeed, a missing trace $B$ on a $(d+1)$-set $T$ gives the relation $x_B\prod_{i\in T\setminus B}(x_i-1)=0$ on $\mathcal A$, expressing $x_T$ in lower degrees; multiplication and induction reduce every higher-degree monomial.

We prove Theorem~\ref{thm:antichain-d2}. The following lemma settles the case $d=1$ and provides the main auxiliary ingredient for $d=2$.

\begin{lemma}\label{lem:app-vc1}
Let $\mathcal A$ be a nonempty antichain of nonempty sets with $\operatorname{VCdim}(\mathcal A)\le1$, and put $U=\operatorname{supp}\mathcal A$. Then $|\mathcal A|\le|U|$. If every member has size at least two, equality holds if and only if $\mathcal A=\binom U{|U|-1}$ with $|U|\ge3$.
\end{lemma}

\begin{proof}
The case $|\mathcal A|=1$ is immediate. Fix $H\in\mathcal A$, and let $\delta_H$ denote its point indicator on $\mathcal A$. Since $\mathcal A$ is an antichain, the projection onto $U\setminus H$ has the empty trace only at $H$. Thus the indicator of the empty set on the projected family pulls back to $\delta_H$. As the projected family has VC-dimension at most one, VC-one interpolation gives $\delta_H=1+\sum_{i\notin H}c_{H,i}x_i$ on $\mathcal A$. The functions $1-\delta_H$ span $\mathbb R^{\mathcal A}$, since their evaluation matrix is $J-I$. Thus the coordinate functions $x_i$, $i\in U$, span this space, proving the bound.

Suppose equality holds and every member has size at least two. The coordinate functions form a basis. Write $1=\sum_{i\in U}h_ix_i$,
and let $m_i$ and $m_{ij}$ count members containing $i$ and $\{i,j\}$, respectively. Then $\delta_H=1+\sum_{i\notin H}c_{H,i}x_i=\sum_{i\in U}h_ix_i+\sum_{i\notin H}c_{H,i}x_i=\sum_{i\in H}h_ix_i+\sum_{i\notin H}(c_{H,i}+h_i)x_i$. For $i\in H$, the coefficient of $x_i$ in $\delta_H$ is $h_i$. Reading the $x_i$-coordinate in $x_i=\sum_{H\ni i}\delta_H$ gives $h_i=1/m_i>0$. Hence, whenever $\{i,j\}\in\partial_2\mathcal A$, the $x_i$- and $x_j$-coefficients of $x_ix_j$ are $m_{ij}/m_i>0$ and $m_{ij}/m_j>0$. Neither singleton trace on $\{i,j\}$ can be missing, since its absence would give $x_ix_j=x_i$ or $x_ix_j=x_j$. The full trace occurs, so the VC bound forces the empty trace of $\{i,j\}\in\partial_2\mathcal A$ to be missing. Thus distinct members differ by exactly one point in each direction and have a common size $s\ge2$.

Fix two members $S\cup\{a\}$ and $S\cup\{b\}$, where $|S|=s-1$. Every other member has the form $S\cup\{c\}$ with $c\notin S\cup\{a,b\}$, or
$(S\setminus\{z\})\cup\{a,b\}$ with $z\in S$. These two types cannot coexist, since they differ by two points in each direction. Hence either all members contain $S$, giving $|\mathcal A|=|U|-s+1<|U|$, or all lie in $\binom U{s}$ with $U=S\cup\{a,b\}$. Equality forces the latter family to be the entire layer. The converse follows directly.
\end{proof}

\begin{lemma}\label{lem:app-local-deficit}
Let $\mathcal P\subseteq2^R$ be nonempty, with $R=\operatorname{supp}\mathcal P$ and $\operatorname{VCdim}(\mathcal P)\le2$, and define
$\alpha(\mathcal P)=1+|R|+|\partial_2\mathcal P|-|\mathcal P|$. Suppose $B\subseteq R$ and, for every $z\in B$, the members of $\mathcal P$ containing $z$ form an antichain and have size at least two. Then $\alpha(\mathcal P)\ge|B|$.
\end{lemma}

\begin{proof}
Induct on $|B|$. When $B=\varnothing$, VC interpolation gives $|\mathcal P|\le1+|R|+|\partial_2\mathcal P|$. Otherwise choose $z\in B$, and put $\mathcal P_0=\{P\in\mathcal P:z\notin P\}$, $\mathcal P_1=\{P\setminus\{z\}:P\in\mathcal P,\ z\in P\}$, $\mathcal P'=\mathcal P_0\cup\mathcal P_1$, and $\mathcal D=\mathcal P_0\cap\mathcal P_1$. Let $N(z)=\{u\in R\setminus\{z\}:\{z,u\}\in\partial_2\mathcal P\}$.
The support of $\mathcal P'$ is $R\setminus\{z\}$. Counting members and shadow pairs gives
\[
\alpha(\mathcal P)=\alpha(\mathcal P')+1+|N(z)|-|\mathcal D|.
\]
The family $\mathcal D\subseteq\mathcal P_1$ is an antichain of nonempty sets. Its VC-dimension is at most one, since a pair shattered by $\mathcal D$ would extend to a shattered triple in $\mathcal P$ by adjoining $z$. Moreover, $\operatorname{supp}\mathcal D\subseteq N(z)\setminus B$: if $D\in\mathcal D$ contained $w\in B$, then $D,D\cup\{z\}\in\mathcal P$ would be comparable members containing $w$.
Lemma~\ref{lem:app-vc1} gives $|\mathcal D|\le|N(z)\setminus B|$, also when $\mathcal D=\varnothing$.

Set $B'=B\setminus(\{z\}\cup N(z))$. For $w\in B'$, we have $w\notin N(z)$, so no member of $\mathcal P$ contains both $w$ and $z$. Thus deleting $z$ leaves the subfamily containing $w$ unchanged. Therefore, $\mathcal P'$ satisfies the hypotheses with marked set $B'$.
Induction yields
\[
\alpha(\mathcal P) \ge\alpha(\mathcal P')+1+|N(z)\cap B| \ge|B'|+1+|N(z)\cap B|=|B|.
\]
\end{proof}

\begin{proof}[Proof of Theorem~\ref{thm:antichain-d2}]
For $d=1$, \cref{lem:app-vc1} gives $|\mathcal A|\le|\operatorname{supp}\mathcal A| =|\partial_1\mathcal A|$, with equality precisely when
$\mathcal A=\binom Qr$ for $|Q|=r+1$ and $r\ge2$. This is the required assertion for $d=1$.

We henceforth assume $d=2$, write $\mathcal F=\mathcal A$, and put $Q=\operatorname{supp}\mathcal F$. Let $G$ be the graph on $Q$ with edge set $\partial_2\mathcal F$. Fix an arbitrary edge $\{a,b\}$ and put $R=Q\setminus\{a,b\}$. For $\varepsilon,\eta\in\{0,1\}$, let $\mathcal F_{\varepsilon\eta}$ consist of the sets $F\setminus\{a,b\}$ with $F\in\mathcal F$, $\mathbf1_{\{a\in F\}}=\varepsilon$, and $\mathbf1_{\{b\in F\}}=\eta$. Set $\mathcal P=\bigcup_{\varepsilon,\eta}\mathcal F_{\varepsilon\eta}$ and $\mathcal D=\mathcal F_{10}\cap\mathcal F_{01}$.
Since $\mathcal F$ is an antichain, the only possible overlap between distinct sections is $\mathcal F_{10}\cap\mathcal F_{01}$. Hence $|\mathcal F|=|\mathcal P|+|\mathcal D|$. Also $\operatorname{VCdim}(\mathcal P)\le2$, while $\mathcal D$ is an antichain of sets of size at least two with VC-dimension at most one: a shattered pair in $\mathcal D$, together with $a$, would be shattered by $\mathcal F$.

Let $V_a=N_G(a)\setminus\{b\}$, $V_b=N_G(b)\setminus\{a\}$, $K=V_a\cap V_b$, and $B=R\setminus(V_a\cup V_b)$. Then every original member containing $z\in B$ avoids $a,b$. Thus the members of $\mathcal P$ containing $z$ form an antichain and retain their original sizes, all at least three. Lemma~\ref{lem:app-local-deficit} gives $\alpha(\mathcal P)\ge|B|$. Since $\operatorname{supp}\mathcal D\subseteq K$, Lemma~\ref{lem:app-vc1} also gives $|\mathcal D|\le|K|$. Counting shadow pairs according to their intersections with $\{a,b\}$, we obtain
\[
\begin{aligned}
|\partial_2\mathcal F|-|\mathcal F|
&=|\partial_2\mathcal P|+|V_a|+|V_b|+1-|\mathcal P|-|\mathcal D|\\
&=\bigl(\alpha(\mathcal P)-|B|\bigr)+\bigl(|K|-|\mathcal D|\bigr)\ge0.
\end{aligned}
\]
This proves the inequality.

Suppose equality holds. Since the edge $\{a,b\}$ was arbitrary, $|\mathcal D|=|K|$ for every edge. Each edge belongs to a member of size at least three, so $K=N_G(a)\cap N_G(b)\ne\varnothing$. Equality in $|\mathcal D|\le|\operatorname{supp}\mathcal D|\le|K|$ and Lemma~\ref{lem:app-vc1} imply $\mathcal D=\binom K{|K|-1}$ and $|K|\ge3$. Every pair in $K$ lies in some $D\in\mathcal D$; as $D\cup\{a\}\in\mathcal F$, this pair is an edge of $G$. Thus the common neighborhood of every edge of $G$ is a clique.

Choose a maximal clique $C$ of $G$ containing an edge. Let $a,b$ be any two distinct vertices of $C$, and put $K=N_G(a)\cap N_G(b)$. Since $C$ is a clique, we have $C\setminus\{a,b\}\subseteq K$. Moreover, $K$ is a clique by the preceding argument, so $\{a,b\}\cup K$ is a clique containing $C$. By maximality, $K=C\setminus\{a,b\}$. Set $r=|K|=|C|-2\ge3$.

We claim that $\binom Cr =\mathcal F$. First, let $H\in\binom Cr$ and choose $a\in H$ and $b\in C\setminus H$. For this choice of $a,b$, the preceding conclusion gives $K=N_G(a)\cap N_G(b)=C\setminus\{a,b\}$. Hence $H\setminus\{a\}\subseteq K$ and $|H\setminus\{a\}|=r-1$. By the equality characterization established above for the edge $\{a,b\}$, every $(r-1)$-subset $D$ of $K$ satisfies $D\cup\{a\}\in\mathcal F$. Taking $D=H\setminus\{a\}$ yields $H\in\mathcal F$. Since $H$ was arbitrary, we conclude that $\binom Cr\subseteq\mathcal F$.

It remains to prove the reverse inclusion. Suppose that $F\in\mathcal F\setminus\binom Cr$. Since $\mathcal F$ is an antichain and $\binom Cr\subseteq\mathcal F$, we have $|F\cap C|\le r-1$, and hence $|C\setminus F|\ge3$. Choose $T\in\binom{C\setminus F}{3}$. Since $|C\setminus T|=r-1$, the complete layer $\binom Cr$ realizes every nonempty trace on $T$. Together with $F\cap T=\varnothing$, this implies that $T$ is shattered by $\mathcal F$, contradicting $\operatorname{VCdim}(\mathcal F)\le2$. Therefore $\mathcal F=\binom Cr$.
Conversely, every complete layer $\binom Cr$ with $|C|=r+2$ and $r\ge3$ has VC-dimension two and satisfies $\partial_2\mathcal F=\binom C2$. Hence $|\mathcal F|=\binom{r+2}r=\binom{r+2}2=|\partial_2\mathcal F|$, completing the proof.
\end{proof}
\end{appendices}

\end{document}